\documentclass[12pt]{amsart}

\usepackage{graphicx}
\usepackage{amsmath}

\def\al{\alpha}

\def\ga{\gamma}

\def\la{\lambda}

\def\om{\omega}

\newcommand{\der}{{\rm d}}

\numberwithin{equation}{section}

\newtheorem{theorem}{Theorem}[section]

\newtheorem{corollary}{Corollary}[theorem]
\theoremstyle{remark}

\theoremstyle{remark}

\usepackage{amssymb,stmaryrd}
\usepackage{amscd}

\author{Matthew Randall}

\email{mran052@gmail.com}

\title{Automorphisms and transformations of solutions to the generalised Chazy equation for various parameters}


\begin{document}

\begin{abstract}
We analyse the automorphisms of solutions to Chazy's equation and the generalised Chazy equation for the parameters $k=2,3,4$ and $9$. These automorphisms are induced by triangular domains with isosceles symmetry. We also prove theorems about the transformations of solutions to the generalised Chazy equation between various different parameters. Using the transformation of solutions between parameters $k=2$ and $k=\frac{2}{3}$, we are able to prove a result about the automorphism of the solutions to the $k=\frac{2}{3}$ generalised Chazy equation. 
\end{abstract}

\maketitle

\pagestyle{myheadings}
\markboth{Randall}{Automorphisms and transformations of solutions to the generalised Chazy equation}

\section{Introduction}
In the earlier paper \cite{r18} we studied automorphisms of solutions to Ramanujan's differential equations arising from the equilateral symmetry of the domain of the Schwarz triangle function $s(0,0,0,x)$. We also studied the automorphisms of the solutions to the $k=\frac{3}{2}$ generalised Chazy equation, induced by the equilateral spherical triangle with angles $(\frac{2}{3}\pi,\frac{2}{3}\pi,\frac{2}{3}\pi)$, and showed its relationship to solutions of the $k=3$ generalised Chazy equation. We now continue to investigate the situation when the automorphisms of solutions to Chazy's equation and the generalised Chazy equations arise from triangles with isosceles symmetry. This happens in the generalised Chazy equation when the parameters are $k=2,3,4$ and $9$. In addition, the occurrence of the same Schwarz triangle function in the solution to the generalised Chazy equation with various different parameters gives rise to transformations from the solution of one parameter to another. In such transformations, the domains of the Schwarz triangle functions that occur in the solutions do not necessarily need to have isosceles symmetry. 

We study the automorphisms of solutions to Chazy's equations once again in Section \ref{autinf}.
In Section \ref{aut2}, we give automorphisms of solutions to $k=2$ generalised Chazy equation and deduce the transformations that map these solutions to solutions of the equations with parameters $k=3$ and $k=\frac{2}{3}$. Furthermore, we use the transformation of solutions between the parameters $k=2$ and $k=\frac{2}{3}$ to deduce automorphisms of solutions to the $k=\frac{2}{3}$ generalised Chazy equation.

In Section \ref{aut3}, we deduce the automorphisms of solutions to the $k=3$ generalised Chazy equation and the transformations that map the solutions to the $k=3$ generalised Chazy equation to the solutions with parameters $k=2,4$ and $k=\frac{3}{2}$ (the transformation to the solution of the equation with parameter $k=\frac{3}{2}$ was already given in \cite{r18}). In Section \ref{aut4}, we give the automorphism of solutions to the $k=4$ generalised Chazy equation and its transformation to solutions of the $k=3$ generalised Chazy equation. We also deduce the automorphisms of the solutions to the $k=9$ generalised Chazy equation and the transformation to the one with parameter $k=18$ in Section \ref{aut9}. 

The generalised Chazy equation with parameter $k$ is given by
\begin{align*}
y'''-2yy''+3y'^2-\frac{4}{36-k^2}(6y'-y^2)^2=0
 \end{align*} 
and Chazy's equation 
\begin{align*}
y'''-2yy''+3y'^2=0
 \end{align*} 
is obtained in the limit as $k$ tends to infinity. The generalised Chazy equation was introduced in \cite{chazy1}, \cite{chazy2} and studied more recently in \cite{co96}, \cite{acht}, \cite{ach} and \cite{bc17}. The generalised Chazy equation with parameters $k=\frac{2}{3}$, $\frac{3}{2}$, $2$ and $3$ was also further investigated in \cite{r16}. The solution to the generalised Chazy equation is given by the following. Let 
\begin{align}\label{weq}
w_1&=-\frac{1}{2}\frac{\der }{\der x}\frac{s'}{s(s-1)},\nonumber \\
w_2&=-\frac{1}{2}\frac{\der }{\der x}\frac{s'}{s-1},\\
w_3&=-\frac{1}{2}\frac{\der }{\der x}\frac{s'}{s},\nonumber
\end{align}
where $s=s(\al,\beta,\ga,x)$ is a solution to the Schwarzian differential equation 
\[
\{s,x\}+\frac{1}{2}(s')^2V=0
\]
and
\[
\{s,x\}=\frac{\der}{\der x}\left(\frac{s''}{s'}\right)-\frac{1}{2}\left(\frac{s''}{s'}\right)^2
\]
is the Schwarzian derivative with the potential $V$ given by
\[
V=\frac{1-\beta^2}{s^2}+\frac{1-\ga^2}{(s-1)^2}+\frac{\beta^2+\ga^2-\al^2-1}{s(s-1)}.
\]
The combination $y=-2w_1-2w_2-2w_3$ solves the generalised Chazy equation when 
\[
(\al,\beta,\ga)=\left(\frac{1}{3},\frac{1}{3},\frac{2}{k}\right) \text{~or~} \left(\frac{2}{k},\frac{2}{k},\frac{2}{k}\right).
\]

The combination $y=-w_1-2w_2-3w_3$ solves the generalised Chazy equation when 
\[
(\al,\beta,\ga)=\left(\frac{1}{k},\frac{1}{3},\frac{1}{2}\right) \text{~or~} \left(\frac{1}{k},\frac{2}{k},\frac{1}{2}\right) \text{~or~}  \left(\frac{1}{k},\frac{1}{3},\frac{3}{k}\right) ,
\]
with permutations of $w_1$, $w_2$ and $w_3$ in $y$ corresponding to permutations of the values $\al$, $\beta$ and $\ga$ in $(\al,\beta,\ga)$.
The combination $y=-w_1-w_2-4w_3$ solves the generalised Chazy equation whenever 
\[
(\al,\beta,\ga)=\left(\frac{1}{k},\frac{1}{k},\frac{4}{k}\right) \text{~or~} \left(\frac{1}{k},\frac{1}{k},\frac{2}{3}\right),
\]
again permutating $w_1$, $w_2$ and $w_3$ in $y$ corresponds to permutating the values $\al$, $\beta$, $\ga$ in $(\al,\beta,\ga)$.
Following \cite{acht}, the functions $w_1$, $w_2$ and $w_3$ satisfy the following system of differential equations:
\begin{align}\label{wde}
w_1'=w_2w_3-w_1(w_2+w_3)+\tau^2,\nonumber\\
w_2'=w_3w_1-w_2(w_3+w_1)+\tau^2,\\
w_3'=w_1w_2-w_3(w_1+w_2)+\tau^2,\nonumber
\end{align}
where
\[
\tau^2=\al^2(w_1-w_2)(w_3-w_1)+\beta^2(w_2-w_3)(w_1-w_2)+\ga^2(w_3-w_1)(w_2-w_3).
\]

Taking $P=y$, the generalised Chazy equation with parameter $k$ is equivalent to the following non-linear first order system of differential equations
\begin{align}\label{ndek0}
\frac{\der P}{\der x}&=\frac{1}{6}(P^2-Q),\nonumber\\
\frac{\der Q}{\der x}&=\frac{2}{3}(PQ-R),\\
\frac{\der R}{\der x}&=PR+\frac{k^2}{36-k^2}Q^2.\nonumber
\end{align}
Consider the triple $(P,Q,R)$ of functions of the variable $x$. We say that $(P,Q,R)$ satisfies or solves the non-linear first order differential equation associated to the generalised Chazy equation with parameter $k$ if $P$, $Q$, $R$ satisfy the system of differential equations (\ref{ndek0}). In this case, $y=P$ is a solution to the generalised Chazy equation with parameter $k$.
For Chazy's equation, we obtain Ramanujan's differential equations given by
\begin{align}\label{rde0}
\frac{\der P}{\der x}&=\frac{1}{6}(P^2-Q),\nonumber\\
\frac{\der Q}{\der x}&=\frac{2}{3}(PQ-R),\\
\frac{\der R}{\der x}&=PR-Q^2.\nonumber
\end{align}
We say that $(P,Q,R)$ satisfies or solves Ramanujan's differential equation associated to Chazy's equation if $P$, $Q$, $R$ satisfy the system of differential equations (\ref{rde0}). In this case, $y=P$ is a solution to Chazy's equation.

\section{Automorphisms of solutions to Ramanujan's differential equations}\label{autinf}

The results in this section are obtained by considering solutions to Chazy's equation associated to Schwarz triangle functions that have triangular domains with isosceles symmetry. We analyse the solutions given by the Schwarz triangle functions $s(0,0,\frac{1}{2},x)$ and $s(0,0,\frac{1}{3},x)$. The analysis of the first function recovers partly the result from the previous paper \cite{r18}. 
We can state it in the following manner:
\begin{theorem}\label{thm0}
Suppose the triple $(P,Q,R)$ satisfies Ramanujan's differential equations (\ref{rde0}). Let  $\nu$ be a solution to the cubic equation
\begin{align*}
\nu^3-\frac{3}{4}Q\nu+\frac{1}{4}R=0. 
 \end{align*}
Then the triple
\begin{align*}
\left(P-\nu,-Q+5\nu^2,\frac{11}{4}R-\frac{21}{4}Q\nu\right)
\end{align*}
also satisfies Ramanujan's differential equations.
\end{theorem}

\begin{proof}
When $s=s(0,0,\frac{1}{2},x)$, taking $y=P=-w_1-2w_2-3w_3$ as a solution to Chazy's equation, we obtain 
$Q=\frac{5}{2}(w_1-w_2)^2-\frac{3}{2}(w_1-w_2)(w_1+w_2-2w_3)$ and $R=\frac{7}{2}(w_1-w_2)^3-\frac{9}{2}(w_1-w_2)^2(w_1+w_2-2w_3)$. We have used equations (\ref{wde}) and the formulas $Q=P^2-6P'$, $R=PQ-\frac{3}{2}Q'$ to determine them with $(\al,\beta,\ga)=(0,0,\frac{1}{2})$. When we take $y=p=-2w_1-w_2-3w_3$, also as a solution to Chazy's equation for the same Schwarz triangle function, we obtain $q=\frac{5}{2}(w_1-w_2)^2+\frac{3}{2}(w_1-w_2)(w_1+w_2-2w_3)$ and $r=-\frac{7}{2}(w_1-w_2)^3-\frac{9}{2}(w_1-w_2)^2(w_1+w_2-2w_3)$ with $q=p^2-6p'$ and $r=pq-\frac{3}{2}q'$.

The next step is to solve for $p$, $q$, $r$ in terms of $P$, $Q$, $R$.
Let $\nu=w_1-w_2$, $\la=w_2-w_3$ and $\mu=w_3-w_1$. We find 
\begin{align*}
Q=&\frac{5}{2}\nu^2-\frac{3}{2}\nu(\la-\mu),\\
R=&\frac{7}{2}\nu^3-\frac{9}{2}\nu^2(\la-\mu).
\end{align*}
Eliminating $\la-\mu$ from these two equations gives us the cubic equation
\begin{equation}\label{cubic1}
\nu^3-\frac{3}{4}Q\nu+\frac{1}{4}R=0.
\end{equation}
Consequently, we obtain the triple
\begin{align*}
(p,q,r)&=\left(P-\nu,\frac{5}{2}\nu^2+\frac{3}{2}\nu(\la-\mu),-\frac{7}{2}\nu^3-\frac{9}{2}\nu^2(\la-\mu)\right)\\
&=\left(P-\nu,-Q+5\nu^2,\frac{11}{4}R-\frac{21}{4}Q\nu\right)
\end{align*}
that also solves Ramanujan's differential equations.
\end{proof}
Since a root of the cubic equation (\ref{cubic1}) is given by
\[
\nu=-\frac{1}{2}(R+\sqrt{R^2-Q^3})^{\frac{1}{3}}-\frac{1}{2}\frac{Q}{(R+\sqrt{R^2-Q^3})^{\frac{1}{3}}},
\]
this recovers part of the result of Theorem 1.1 in \cite{r18}.

The analysis of the solution using the Schwarz triangle function $s(0,0,\frac{1}{3},x)$ gives us a new result. 
We prove the following theorem:
\begin{theorem}\label{thm1}
Suppose the triple $(P,Q,R)$ satisfies Ramanujan's differential equations (\ref{rde0}).
Let  $\nu$ be a solution to the quartic equation
\begin{align*}
\nu^4-\frac{2}{3}Q\nu^2+\frac{8}{27}R\nu-\frac{1}{27}Q^2=0. 
 \end{align*}
Then the triple
\begin{align*}
\left(P-2\nu,-Q+10\nu^2,-R-35\nu^3+7Q\nu\right)
 \end{align*}
also satisfies Ramanujan's differential equations.
\end{theorem}

\begin{proof}
When $s=s(0,0,\frac{1}{3},x)$, for the solution of Chazy's equation given by $y=P=-w_1-3w_2-2w_3$, we obtain 
$Q=(w_1-w_2)((w_1-w_2)+8(w_3-w_2))$ and $R=-(w_1-w_2)((w_1-w_3)^2-9(w_2-w_3)^2+18(w_1-w_2)(w_2-w_3))$ with $Q=P^2-6P'$ and $R=PQ-\frac{3}{2}Q'$.

For the solution of Chazy's equation given by $y=p=-3w_1-w_2-2w_3$, we obtain 
$q=(w_1-w_2)(8(w_1-w_3)+(w_1-w_2))$ and $r=-(w_1-w_2)(-(w_2-w_3)^2+9(w_1-w_3)^2+18(w_1-w_3)(w_1-w_2))$ with $q=p^2-6p'$ and $r=pq-\frac{3}{2}q'$.
Once again we solve for $p$, $q$, $r$ in terms of $P$, $Q$, $R$.
Let $\nu=w_1-w_2$, $\la=w_2-w_3$ and $w_3-w_1=-\la-\nu$. Then we get
\begin{align*}
Q=&\nu(\nu-8\la)=\nu^2-8\nu\la,\\
R=&-\nu((\nu+\la)^2-9\la^2+18\nu\la).
\end{align*}
We find through eliminating $\la$ that $\nu$ satisfies the quartic equation
\[
\nu^4-\frac{2}{3}Q\nu^2+\frac{8}{27}R\nu-\frac{1}{27}Q^2=0
\]
with $\la=\frac{1}{8}\nu-\frac{Q}{8\nu}$.
We therefore obtain the triple
\begin{align*}
(p,q,r)&=\left(P-2\nu,\nu(8(\nu+\la)+\nu),-\nu(-\la^2+9(\la+\nu)^2+18\nu(\nu+\la))\right)\\
&=\left(P-2\nu,9\nu^2+8\nu\la,-R-\nu(-10\la^2+10(\la+\nu)^2+18\nu^2+36\nu\la)\right)\\
&=\left(P-2\nu,9\nu^2+8\nu\la,-R-28\nu^3-56\nu^2\la\right)\\
&=\left(P-2\nu,-Q+10\nu^2,-R-35\nu^3+7Q\nu\right)
\end{align*}
solving Ramanujan's differential equations. 
\end{proof}

\section{Automorphisms of solutions to $k=2$ generalised Chazy equation and transformations to solutions of the $k=3$ and $k=\frac{2}{3}$ generalised Chazy equation}\label{aut2}
This section of the paper concerns automorphisms of solutions to the generalised Chazy equation with parameter $k=2$ and the transformations of these solutions to the solutions of the generalised Chazy equation with parameters $k=3$ and $k=\frac{2}{3}$. When $k=2$ we have the following system of equations from (\ref{ndek0}):
\begin{align}\label{nde1}
\frac{\der P}{\der x}&=\frac{1}{6}(P^2-Q),\nonumber\\
\frac{\der Q}{\der x}&=\frac{2}{3}(PQ-R),\\
\frac{\der R}{\der x}&=PR+\frac{1}{8}Q^2.\nonumber
\end{align}
Taking $y=P$ gives a solution to the generalised Chazy equation 
\begin{equation}\label{chazy32}
y'''-2yy''+3y'^2-\frac{1}{8}(6y'-y^2)^2=0
\end{equation}
with parameter $k=2$. We prove the following theorem. 
\begin{theorem}\label{thm2}
Suppose $(P,Q,R)$ satisfies (\ref{nde1}) the non-linear system of differential equations associated to the generalised Chazy equation with parameter $k=2$.
Then the following holds. The triples
\begin{align*}
\left(P+\frac{Q^2}{4R},-Q-\frac{Q^4}{8R^2},-R-\frac{Q^3}{4R}-\frac{Q^6}{64 R^3}\right)
 \end{align*}
and
\begin{align*}
\left(P-2 \nu,\frac{5}{3}Q+\frac{64}{3}\mu^2-8\nu^2,-\frac{7}{27}R-\frac{4}{3}Q\mu-2Q\nu-\frac{Q^2}{R}\nu^2-\frac{2}{27}\frac{Q^3}{R}\right)
\end{align*} 
both satisfy the same system of differential equations (\ref{nde1}), where $\nu$
and $\mu$ are roots of the cubic equations
\begin{align*}
\mu^3+\frac{3}{32}Q\mu+\frac{R}{32}=0,\\
\nu^3-\frac{1}{8}\frac{Q^2}{R}\nu^2-\frac{1}{108}\frac{Q^3}{R}-\frac{2}{27}R=0,
\end{align*}
and $\nu$ is related to $\mu$ by $\nu=-\frac{4}{3}\mu-\frac{Q}{24\mu}$.
Furthermore let 
\[
\chi=\frac{R}{2Q}\pm\frac{1}{2}\sqrt{\frac{R^2}{Q^2}+\frac{Q}{8}}.
\]
Then the triple
\[
\left(P+\chi,\frac{3}{4}Q-3\chi^2,\frac{3}{4}R+\frac{9}{8}Q\chi+3\chi^3\right)
\]
solves the first order system of differential equations associated to the $k=3$ generalised Chazy equation.
\end{theorem}

\begin{proof}
The result is proved by studying the solutions to the generalised Chazy equation given by the Schwarz triangle functions $s(\frac{1}{2},\frac{1}{2},\frac{1}{3},x)$,  $s(\frac{1}{2},\frac{1}{2},1,x)$ and  $s(\frac{1}{3},\frac{1}{3},1,x)$.
For the Schwarz triangle function $s(\frac{1}{2},\frac{1}{2},\frac{1}{3},x)$, both $y=-w_1-3w_2-2w_3$ and $y=-3w_1-w_2-2w_3$ for $w_1$, $w_2$ and $w_3$ given by (\ref{weq}) solve the $k=2$ generalised Chazy equation. 
For $y=P=-w_1-3w_2-2w_3$ we obtain the triple
\[
(P,Q,R)=(-w_1-3w_2-2w_3,-8(w_1-w_3)(w_1-w_2),8(w_1-w_3)^2(w_1-w_2))
\]
from applying (\ref{wde}) with $(\al,\beta,\ga)=(\frac{1}{2},\frac{1}{2},\frac{1}{3})$ and using $Q=P^2-6P'$ and $R=PQ-\frac{3}{2}Q'$. In this instance we obtain
\begin{align}\label{w223}
w_1-w_3&=-\frac{R}{Q},\nonumber\\
w_1-w_2&=\frac{Q^2}{8R},\\
w_2-w_3&=-\frac{Q^2}{8R}-\frac{R}{Q}.\nonumber
\end{align}
For  $y=p=-3w_1-w_2-2w_3$ we have the triple
\[
(p,q,r)=(-3w_1-w_2-2w_3,8(w_2-w_3)(w_1-w_2),-8(w_2-w_3)^2(w_1-w_2))
\]
from again applying (\ref{wde}) and upon substitution of (\ref{w223}) above we obtain
\[
(p,q,r)=\left(P+\frac{Q^2}{4R},-Q-\frac{Q^4}{8R^2},-R-\frac{Q^3}{4R}-\frac{Q^6}{64R^3}\right).
\]

Next, for the Schwarz triangle function $s(\frac{1}{2},\frac{1}{2},1,x)$, both $y=-w_1-3w_2-2w_3$ and $y=-3w_1-w_2-2w_3$ for $w_1$, $w_2$ and $w_3$ given by (\ref{weq}) solve the $k=2$ generalised Chazy equation. For $y=P=-w_1-3w_2-2w_3$, we get the triple
\[
(P,Q,R)=(-w_1-3w_2-2w_3,8(w_3-w_1)(w_1-w_3+3(w_2-w_3)),-8(w_3-w_1)^2(8(w_2-w_3)+(w_2-w_1)))
\]
from (\ref{wde}) while for $y=p=-3w_1-w_2-2w_3$ we have
\[
(p,q,r)=(-3w_1-w_2-2w_3,-8(w_2-w_3)(w_2-w_3+3(w_1-w_3)),8(w_2-w_3)^2(8(w_3-w_1)+(w_2-w_1))).
\]
Now take $\nu=w_1-w_2$ and $\mu=w_3-w_1$, with $w_2-w_3=-\mu-\nu$.
Then we find 
\begin{align*}
Q=&8\mu(-\mu+3(-\nu-\mu))=-32\mu^2-24\mu\nu,\\
R=&-8\mu^2(8(-\mu-\nu)-\nu)=64\mu^3+72\mu^2\nu,
\end{align*}
which give the cubic equation
\begin{align*}
\mu^3+\frac{3}{32}Q\mu+\frac{R}{32}=0
\end{align*}
and  $\nu=-\frac{4}{3}\mu-\frac{Q}{24\mu}$, where $\nu$ satisfies the cubic equation
\begin{align*}
\nu^3-\frac{Q^2}{8R}\nu^2-\frac{1}{108}\frac{Q^3}{R}-\frac{2}{27}R=0.
\end{align*}
We therefore obtain
\begin{align*}
(p,q,r)&=(P-2\nu,-8(\mu+\nu)(4\mu+\nu),8(\mu+\nu)^2(8\mu-\nu))\\
&=(P-2\nu,-32\mu^2-8\nu^2-40\mu\nu,64\mu^3+120\mu^2\nu+48\mu\nu^2-8\nu^3)\\
&=\left(P-2 \nu,\frac{5}{3}Q+\frac{64}{3}\mu^2-8\nu^2,-\frac{7}{27}R-\frac{4}{3}Q\mu-2Q\nu-\frac{Q^2}{R}\nu^2-\frac{2}{27}\frac{Q^3}{R}\right)
\end{align*}
as another solution to (\ref{nde1}). In the last step we have used both cubic equations, the equations for $Q$ and $R$ and the formula
\begin{align*}
\mu\nu^2=\frac{1}{54}(64\mu^3-R-\frac{9}{4}Q\nu).
\end{align*}

Lastly, consider the Schwarz triangle function $s(\frac{1}{3},\frac{1}{3},1,x)$. The function $y=-2w_1-2w_2-2w_3$ solves the $k=2$ generalised Chazy equation while $y=-w_1-2w_2-3w_3$ and $y=-2w_1-w_2-3w_3$ both solve the $k=3$ generalised Chazy equation. For $y=P=-2w_1-2w_2-2w_3$, we obtain the triple
\[
(P,Q,R)=(-2w_1-2w_2-2w_3,-32(w_1-w_3)(w_2-w_3),-32(w_2-w_3)(w_1-w_3)(w_1+w_2-2w_3))
\]
that solves the first order system of differential equations associated to the $k=2$ generalised Chazy equation. 
We also obtain the triples
\begin{align*}
(p,q,r)=&(-w_1-2w_2-3w_3,-3(w_1-w_3)(w_1-w_3+8(w_2-w_3)),\\
&~3(w_1-w_3)((w_1-w_2)^2-9(w_2-w_3)^2+18(w_1-w_3)(w_3-w_2)))
\end{align*}
with $y=p=-w_1-2w_2-3w_3$
and
\begin{align*}
(\tilde p,\tilde q,\tilde r)=&(-2w_1-w_2-3w_3,-3(w_2-w_3)(w_2-w_3+8(w_1-w_3)),\\
&-3(w_2-w_3)(-(w_1-w_2)^2+9(w_1-w_3)^2+18(w_1-w_3)(w_2-w_3)))
\end{align*}
with $y=\tilde p=-2w_1-w_2-3w_3$ that both solve the system of differential equaions associated to the $k=3$ generalised Chazy equation. Let $\la=w_2-w_3$ and $-\mu=w_1-w_3$. Consequently, we obtain
\begin{align*}
\frac{R}{Q}&=\la-\mu,\\
-\la \mu&=-32Q. 
\end{align*}
Solving the quadratic equation $\chi^2-\frac{R}{Q}\chi-32Q=0$ for $\chi=-\mu$ and $\chi=\la$ gives $\chi=\frac{R}{2Q}\pm\frac{1}{2}\sqrt{\frac{R^2}{Q^2}+\frac{Q}{8}}$ and substituting this back into the formulas for $(p,q,r)$ and $(\tilde p, \tilde q, \tilde r)$ above gives the triple
\[
\left(P+\chi,\frac{3}{4}Q-3\chi^2,\frac{3}{4}R+\frac{9}{8}Q\chi+3\chi^3\right)
\]
This gives the desired transformation from $(P,Q,R)$ to the triple solving the first order system associated to the $k=3$ generalised Chazy equation. 
\end{proof}

We now find the transformation of solutions between the generalised Chazy equation with parameters $k=\frac{2}{3}$ and $k=2$.
\begin{theorem}\label{thm23}
Suppose $(P,Q,R)$ satisfies (\ref{nde1}) the non-linear system of differential equations associated to the generalised Chazy equation with parameter $k=2$.
Then the following holds. The triple
\begin{align*}
\left(P+2\mu,\frac{5}{4}Q+10\mu^2,\frac{5}{4}R+\frac{25}{8}Q\mu+35\mu^3\right)
\end{align*} 
satisfies the non-linear system of differential equations associated to the generalised Chazy equation with parameter $k=\frac{2}{3}$, where $\mu$ is a root of the quartic equation
\begin{align*}
\mu^4+\frac{1}{12}Q\mu^2+\frac{1}{27}R\mu-\frac{1}{1728}Q^2=0.
\end{align*}
Conversely, let $(P,Q,R)$ be a solution the non-linear system of differential equations associated to the $k=\frac{2}{3}$ generalised Chazy equation. Then 
\begin{align*}
\left(P-\frac{Q^2}{40R},\frac{4}{5}Q-\frac{Q^4}{800R^2},\frac{4}{5}R-\frac{Q^3}{40R}-\frac{Q^6}{64000R^3}\right)
\end{align*} 
solves the non-linear system of differential equations associated to the $k=2$ generalised Chazy equation.
\end{theorem}

\begin{proof}
We make use of the Schwarz triangle function $s=s(\frac{1}{2},\frac{1}{3},\frac{3}{2},x)$ that appears in both the solutions to the $k=2$ and $k=\frac{2}{3}$ generalised Chazy equation. Note that the domain of this Schwarz triangle function is not isosceles. For the Schwarz triangle function $s(\frac{1}{2},\frac{1}{3},\frac{3}{2},x)$, $y=-w_1-2w_2-3w_3$ solves the $k=2$ generalised Chazy equation while $y=-3w_1-2w_2-w_3$ solves the $k=\frac{2}{3}$ generalised Chazy equation. 
For $y=P=-w_1-2w_2-3w_3$ we obtain the triple
\begin{align*}
&(P,Q,R)=(-w_1-2w_2-3w_3,-8(w_1-w_3)(w_1-w_3+8(w_2-w_3)),\\
&~8(w_1-w_3)((w_1-w_2)^2-9(w_2-w_3)^2-18(w_1-w_3)(w_2-w_3)))
\end{align*}
and for  $y=p=-3w_1-2w_2-w_3$ we obtain the triple
\[
(p,q,r)=(-3w_1-2w_2-w_3,-80(w_1-w_3)(w_2-w_3),-80(w_2-w_3)^2(w_1-w_3)).
\]
Now take $\mu=w_3-w_1$, $\la=w_2-w_3$, and solve for $p$, $q$, $r$ in terms of $P$, $Q$, $R$.
We find
\begin{align*}
Q=&8\mu(-\mu+8\la)=-8\mu^2+64\la \mu,\\
R=&-8\mu((\mu+\la)^2-9\la^2+18\mu\la)=-28 \mu^3-\frac{5}{2}Q\mu+64\mu\la^2.
\end{align*}
Solving this for $\mu$ we get the quartic equation
\[
\mu^4+\frac{1}{12}Q\mu^2+\frac{1}{27}R\mu-\frac{Q^2}{1728}=0
\]
and $\la=\frac{\mu}{8}+\frac{Q}{64\mu}$.
From this we find 
\begin{align*}
p&=P+2\mu,\\
q&=80\mu\la=\frac{5}{4}Q+10\mu^2,\\
r&=80\la^2\mu=\frac{5}{4}R+35\mu^3+\frac{25}{8}Q\mu,
\end{align*}
which gives the triple $(p,q,r)$ a solution to the first order system associated to the $k=\frac{2}{3}$ generalised Chazy equation. 

Conversely, starting from the triple $(P,Q,R)$ a solution to the system of differential equations associated to the $k=\frac{2}{3}$ generalised Chazy equation, we obtain
\begin{align*}
w_2-w_3=&\frac{R}{Q},\\
w_1-w_3=&-\frac{Q^2}{80R},\\
w_1-w_2=&-\frac{Q^2}{80R}-\frac{R}{Q}.
 \end{align*} 
 From this, we get the triple
\begin{align*}
\left(P-\frac{Q^2}{40R},\frac{4}{5}Q-\frac{Q^4}{800R^2},\frac{4}{5}R-\frac{Q^3}{40R}-\frac{Q^6}{64000R^3}\right)
\end{align*}
that satisfes the system of differential equations associated to the $k=2$ generalised Chazy equation.
\end{proof}
As a corollary of Theorem \ref{thm23}, we obtain the following result about the automorphisms of the solution to the $k=\frac{2}{3}$ generalised Chazy equation. 

\begin{corollary}
Suppose $(P,Q,R)$ satisfies the non-linear system of differential equations associated to the generalised Chazy equation with parameter $k=\frac{2}{3}$. Let  $\om=-\frac{1}{2}+\frac{i\sqrt{3}}{2}$ and $\mu$ be one of the following:
\begin{align*}
\mu=&-\frac{1}{240}\frac{Q^2}{R}+\frac{10^{\frac{1}{3}}}{30}\frac{Q}{R^{\frac{1}{3}}}-\frac{10^{\frac{2}{3}}}{15}R^{\frac{1}{3}}\text{~or~}\\
&-\frac{1}{240}\frac{Q^2}{R}+\om\frac{10^{\frac{1}{3}}}{30}\frac{Q}{R^{\frac{1}{3}}}-\om^2\frac{10^{\frac{2}{3}}}{15}R^{\frac{1}{3}}\text{~or~}\\
&-\frac{1}{240}\frac{Q^2}{R}+\om^2\frac{10^{\frac{1}{3}}}{30}\frac{Q}{R^{\frac{1}{3}}}-\om\frac{10^{\frac{2}{3}}}{15}R^{\frac{1}{3}}.
\end{align*}
Then the following holds. The triple
\begin{align*}
\left(P-\frac{1}{40}\frac{Q^2}{R}+2\mu,Q-\frac{Q^4}{640R^2}+10\mu^2,R-\frac{Q^3}{32R}-\frac{1}{51200}\frac{Q^6}{R^3}+\frac{5}{2}Q\mu-\frac{25}{6400}\frac{Q^4}{R^2}\mu+35\mu^3\right)
\end{align*}
satisfies the same non-linear system of differential equations associated to the generalised Chazy equation with parameter $k=\frac{2}{3}$.

\end{corollary}

\begin{proof}
We start with $(P,Q,R)$ a solution to the system of differential equations associated to the $k=\frac{2}{3}$ generalised Chazy equation. From the converse part of Theorem \ref{thm23} we find
\begin{align*}
\left(P-\frac{Q^2}{40R},\frac{4}{5}Q-\frac{Q^4}{800R^2},\frac{4}{5}R-\frac{Q^3}{40R}-\frac{Q^6}{64000R^3}\right)
\end{align*}
satisfies the system of differential equations associated to the $k=2$ generalised Chazy equation. We apply Theorem \ref{thm23} again to solve for $\mu$ in the quartic equation
\[
\mu^4+\left(\frac{1}{15}Q-\frac{Q^4}{9600R^2}\right)\mu^2+\frac{1}{27}\left(\frac{4}{5}R-\frac{Q^3}{40R}-\frac{Q^6}{64000R^3}\right)\mu-\frac{1}{1728}\left(\frac{4}{5}Q-\frac{Q^4}{800R^2}\right)^2=0.
\]
This quartic equation factorises to give
\begin{align*}
\mu=&-\frac{1}{240}\frac{Q^2}{R}+\frac{10^{\frac{1}{3}}}{30}\frac{Q}{R^{\frac{1}{3}}}-\frac{10^{\frac{2}{3}}}{15}R^{\frac{1}{3}},\\
&-\frac{1}{240}\frac{Q^2}{R}+\om\frac{10^{\frac{1}{3}}}{30}\frac{Q}{R^{\frac{1}{3}}}-\om^2\frac{10^{\frac{2}{3}}}{15}R^{\frac{1}{3}},\\
&-\frac{1}{240}\frac{Q^2}{R}+\om^2\frac{10^{\frac{1}{3}}}{30}\frac{Q}{R^{\frac{1}{3}}}-\om\frac{10^{\frac{2}{3}}}{15}R^{\frac{1}{3}} \hspace{12pt} \text{and} \hspace{12pt} \frac{Q^2}{80R}
\end{align*}
as its roots.
Substituting the last root back into the formula in Theorem \ref{thm23} gives us $(P,Q,R)$ once again. For the remaining roots, substituting the formulas in Theorem \ref{thm23} back in gives us 
\begin{align*}
p=&P-\frac{1}{40}\frac{Q^2}{R}+2\mu,\\
q=&Q-\frac{Q^4}{640R^2}+10\mu^2,\\
r=&R-\frac{Q^3}{32R}-\frac{1}{51200}\frac{Q^6}{R^3}+\frac{5}{2}Q\mu-\frac{25}{6400}\frac{Q^4}{R^2}\mu+35\mu^3
\end{align*}
and the triple $(p,q,r)$ gives us the result about the automorphism of solutions to the $k=\frac{2}{3}$ generalised Chazy equation. 

\end{proof}

\section{Automorphisms of solutions to $k=3$ generalised Chazy equation and transformations to solutions of the generalised Chazy equation with parameters $2$, $4$ and $\frac{3}{2}$}\label{aut3}
This section concerns the automorphisms of the solution $y=P$ to the generalised Chazy equation with parameter $k=3$ and the transformations to the solution of the generalised Chazy equation with parameters $k=2$, $k=4$ and $k=\frac{3}{2}$. The transformation to the solution of the $k=\frac{3}{2}$ generalised Chazy equation was already given in Theorem 3.1 of \cite{r18} and we obtain the same result here by analysing two different Schwarz triangle functions. We start off with the following:
\begin{theorem}\label{thm3}
Suppose the triple $(P,Q,R)$ satisfies the non-linear system of differential equations associated to the generalised Chazy equation with parameter $k=3$, i.e. we have
\begin{align*}
\frac{\der}{\der x}P&=\frac{1}{6}(P^2-Q),\nonumber\\
\frac{\der}{\der x}Q&=\frac{2}{3}(PQ-R),\\
\frac{\der}{\der x}R&=PR+\frac{1}{3}Q^2.\nonumber
\end{align*}
Then the following holds. The triples
\begin{align*}
\left(P+\frac{R}{Q},-Q-3\frac{R^2}{Q^2},R+3\frac{R^3}{Q^3}\right)
 \end{align*}
and
\begin{align*}
\left(P-\nu,\frac{5}{4}Q+\frac{27}{4}\mu^2-3\nu^2,\frac{17}{8}R-3\nu^3+\frac{63}{16}Q\mu+\frac{243}{16}\mu^3\right)
\end{align*} 
both satisfy the same non-linear system of differential equations associated to the generalised Chazy equation with parameter $k=3$, where $\mu$ is a root of the quartic equation
\begin{align*}
 \mu^4+\frac{2}{9} Q\mu^2+\frac{8}{81}R\mu-\frac{Q^2}{243}=0,
\end{align*}
and  $\nu=-\frac{9}{8}\mu-\frac{1}{24\mu}Q$ is a solution to
\[
\nu^4+\frac{R}{Q}\nu^3-\frac{3}{64}Q^2-\frac{9}{64}\frac{R^2}{Q}=0.
\]
For the same $\mu$, the triple 
\[
\left(P+\mu,\frac{4}{3}Q+4\mu^2,\frac{4}{3}R+2Q\mu+10\mu^3\right)
\]
satisfies the non-linear system of differential equations associated to the generalised Chazy equation with parameter $k=2$, while the triple 
\[
\left(P-\frac{Q^2}{3R},-\frac{5}{3}Q-\frac{5}{9}\frac{Q^4}{R^2},-\frac{5}{3}\frac{Q^3}{R}-\frac{5}{3}R-\frac{10}{27}\frac{Q^6}{R^3}\right)
\]
satisfies the non-linear system of differential equations associated to the generalised Chazy equation with parameter $k=4$.
\end{theorem}

\begin{proof}
The theorem is proved by considering the solutions to the generalised Chazy equation given by the Schwarz triangle functions  $s(\frac{2}{3},\frac{1}{3},\frac{1}{3},x)$,  $s(\frac{1}{3},\frac{1}{3},1,x)$ and  $s(\frac{1}{3},\frac{1}{3},\frac{1}{2},x)$. For the Schwarz triangle function $s=s(\frac{2}{3},\frac{1}{3},\frac{1}{3},x)$, $y=P=-4w_1-w_2-w_3$ and $y=p=-2w_1-2w_2-2w_3$ both solve the $k=3$ generalised Chazy equation. Using (\ref{wde}) with $(\al,\beta,\ga)=(\frac{2}{3},\frac{1}{3},\frac{1}{3})$ we obtain the triples
\[
(P,Q,R)=(-4w_1-w_2-w_3,-3(w_2-w_3)^2,-3(w_2-w_3)^2(2 w_1-w_2-w_3))
\]
and
\[
(p,q,r)=(-2w_1-2w_2-2w_3,-12(w_1-w_3)(w_1-w_2),12(w_1-w_3)(w_1-w_2)(2w_1-w_2-w_3)).
\]
We solve for $w_1$, $w_2$ and $w_3$ in the first triple to find
\begin{align*}
w_1=&-\frac{1}{6}\left(P-\frac{R}{Q}\right),\\
w_2=&-\frac{1}{6}\left(P+2\frac{R}{Q}\right)\pm\frac{1}{2}\sqrt{-\frac{Q}{3}},\\
w_3=&-\frac{1}{6}\left(P+2\frac{R}{Q}\right)\mp\frac{1}{2}\sqrt{-\frac{Q}{3}}.
\end{align*}
This gives
\[
(p,q,r)=\left(P+\frac{R}{Q},-Q-3\frac{R^2}{Q^2},R+3\frac{R^3}{Q^3}\right)
\]
as also a solution to the first order system of differential equations associated to the $k=3$ generalised Chazy equation. 

Now once again consider the Schwarz triangle function $s(\frac{1}{3},\frac{1}{3},1,x)$. In Section \ref{aut2} earlier we expressed solutions to the $k=3$ generalised Chazy equation in terms of the solution to the $k=2$ generalised Chazy equation. Here we start from a solution to the $k=3$ generalised Chazy equation and find an automorphism that takes it to another solution of the $k=3$ generalised Chazy equation and also a transformation to a solution of the $k=2$ generalised Chazy equation.  

For the Schwarz triangle function $s=s(\frac{1}{3},\frac{1}{3},1,x)$, both $y=-w_1-2w_2-3w_3$ and $y=-2w_1-w_2-3w_3$ solve the $k=3$ generalised Chazy equation while the function $y=-2w_1-2w_2-2w_3$ solves the $k=2$ generalised Chazy equation. Let $P=-w_1-2w_2-3w_3$ and $p=-2w_1-w_2-3w_3$. We find the triples
\begin{align*}
(P,Q,R)=&(-w_1-2w_2-3w_3,-3(w_1-w_3)(w_1-w_3+8(w_2-w_3)),\\
&~3(w_1-w_3)((w_1-w_2)^2-9(w_2-w_3)^2+18(w_1-w_3)(w_3-w_2)))
\end{align*}
and
\begin{align*}
(p,q,r)=&(-2w_1-w_2-3w_3,-3(w_2-w_3)(w_2-w_3+8(w_1-w_3)),\\
&-3(w_2-w_3)(-(w_1-w_2)^2+9(w_1-w_3)^2+18(w_1-w_3)(w_2-w_3))).
\end{align*}
Let $\nu=w_1-w_2$, $-\mu=w_1-w_3$ and $w_2-w_3=-\mu-\nu$. We obtain
\begin{align*}
Q=&3\mu(-\mu+8(-\mu-\nu))=-27\mu^2-24\mu\nu,\\
R=&-3\mu(\nu^2-9(\mu+\nu)^2-18\mu(\mu+\nu))=24\mu\nu^2+108\mu^2\nu+81\mu^3.
\end{align*}
From these formulas for $Q$ and $R$ we get
\begin{align*}
\mu^2\nu&=\frac{1}{81}(Q\nu+R)-\mu^3,\\
\mu\nu^2&=\frac{3}{64}(Q\mu-\frac{8}{9}Q\nu+27\mu^3)
\end{align*}
and also
\[
Q\nu=-\frac{81}{8}\mu^3-\frac{27}{8}Q\mu-R.
\]
Eliminating $\nu$ from the above two equations for $Q$ and $R$, we find that $\mu$ is a root of the quartic equation
\begin{align*}
\mu^4+\frac{2}{9} Q\mu^2+\frac{8}{81}R\mu-\frac{Q^2}{243}=0,
\end{align*}
and $\nu$ is a solution to
\[
\nu^4+\frac{R}{Q}\nu^3-\frac{3}{64}Q^2-\frac{9}{64}\frac{R^2}{Q}=0,
\]
where $\nu$ is related to $\mu$ by $\nu=-\frac{9}{8}\mu-\frac{1}{24\mu}Q$.
For $p=-2w_1-w_2-3w_3$, we therefore have
\begin{align*}
p&=P-\nu,\\
q&=-3(\nu^2+10\mu\nu+9\mu^2)=\frac{5}{4}Q-3\nu^2+\frac{27}{4}\mu^2,\\
r&=3(-\nu^3+45\mu^2\nu+17\mu\nu^2+27\mu^3)\\
&=-3\nu^3-\frac{11}{24}Q\nu+\frac{5}{3}R+\frac{675}{64}\mu^3+\frac{153}{64}Q\mu\\
&=\frac{17}{8}R-3\nu^3+\frac{63}{16}Q\mu+\frac{243}{16}\mu^3.
\end{align*}

Consequently, we find the triple with $p=-2w_1-w_2-3w_3$ given by
\begin{align*}
(p,q,r)=&\left(P-\nu,\frac{5}{4}Q+\frac{27}{4}\mu^2-3\nu^2,\frac{17}{8}R-3\nu^3+\frac{63}{16}Q\mu+\frac{243}{16}\mu^3\right)
\end{align*}
again solves the first order system associated to the $k=3$ generalised Chazy equation 
and the triple with $\tilde p=-2w_1-2w_2-2w_3$ given by
\begin{align*}
(\tilde p,\tilde q, \tilde r)=&(-2w_1-2w_2-2w_3,-32(w_1-w_3)(w_2-w_3),-32(w_2-w_3)(w_1-w_3)(w_1+w_2-2w_3))\\
=&(P+\mu,-32\mu(\mu+\nu),32\mu(\mu+\nu)(2\mu+\nu))\\
=&\left(P+\mu,\frac{4}{3}Q+4\mu^2,\frac{4}{3}R+2Q\mu+10\mu^3\right)
\end{align*}
solves the first order system associated to the $k=2$ generalised Chazy equation. 

Finally to obtain the transformation to solutions of the $k=4$ generalised Chazy equation, consider the Schwarz triangle function $s(\frac{1}{3},\frac{1}{3},\frac{1}{2},x)$. The combinations $y=-w_1-2w_2-3w_3$ and $y=-2w_1-w_2-3w_3$ both solve the $k=3$ generalised Chazy equation while the function $y=-2w_1-2w_2-2w_3$ solves the $k=4$ generalised Chazy equation. Let $P=y=-w_1-2w_2-3w_3$ and $p=y=-2w_1-w_2-3w_3$. We find
\[
(P,Q,R)=(-w_1-2w_2-3w_3,-3(w_1-w_3)(w_1-w_2),3(w_1-w_2)^2(w_1-w_3))
\]
and
\[
(p,q,r)=(-2w_1-w_2-3w_3,3(w_2-w_3)(w_1-w_2),3(w_1-w_2)^2(w_2-w_3)).
\]
Also, the triple 
\[
(\tilde p,\tilde q,\tilde r)=(-2w_1-2w_2-2w_3,-5(w_1-w_3)(w_2-w_3),-5(w_1-w_3)(w_2-w_3)(w_1+w_2-2w_3))
\]
with $\tilde p=-2w_1-2w_2-2w_3$
satisfies the first order system associated to the $k=4$ generalised Chazy equation. 
From the formulas for $Q$ and $R$ we get 
\begin{align*}
w_1-w_2&=-\frac{R}{Q},\\
w_1-w_3&=-\frac{Q^2}{3R},\\
w_2-w_3&=\frac{R}{Q}+\frac{Q^2}{3R}.
\end{align*}
Consequently, we find
\[
(p,q,r)=\left(P+\frac{R}{Q},-Q-3\frac{R^2}{Q^2},R+3\frac{R^3}{Q^3}\right)
\]
solves the $k=3$ equation which we have already deduced earlier from the Schwarz triangle function $s(\frac{2}{3},\frac{1}{3},\frac{1}{3},x)$
and
\[
(\tilde p,\tilde q, \tilde r)=\left(P-\frac{Q^2}{3R},-\frac{5}{3}Q-\frac{5}{9}\frac{Q^4}{R^2},-\frac{5}{3}\frac{Q^3}{R}-\frac{5}{3}R-\frac{10}{27}\frac{Q^6}{R^3}\right)
\]
solves the $k=4$ generalised Chazy equation. 
\end{proof}

This part of the section concerns the transformation of the solution to the generalised Chazy equation with parameter $k=3$ to the solution with parameter $k=\frac{3}{2}$. We recover the result of Theorem 3.1 from \cite{r18}.
\begin{theorem}[Theorem 3.1 of \cite{r18}]\label{thm3}
Suppose $(P,Q,R)$ satisfies the non-linear system of differential equations associated to the generalised Chazy equation with parameter $k=3$.
Then the following holds. The triple
\begin{align*}
\left(P-w,\frac{5}{3}Q+5w^2,\frac{5}{12}R+\frac{5}{12}Qw\right)
\end{align*} 
satisfies the non-linear system of differential equations associated to the generalised Chazy equation with parameter $k=\frac{3}{2}$, where $w$ is a root of the cubic equation
\begin{align*}
w^3+\frac{1}{4}Qw-\frac{1}{12}R=0.
\end{align*}
Conversely, let $(P,Q,R)$ be a solution to the non-linear system of differential equations associated to the $k=\frac{3}{2}$ generalised Chazy equation. Then 
\begin{align*}
\left(P+\frac{R}{Q},\frac{3}{5}Q-3\frac{R^2}{Q^2},\frac{9}{5}R+3\frac{R^3}{Q^3}\right).
\end{align*}
solves the non-linear system of differential equations associated to the $k=3$ generalised Chazy equation.
\end{theorem}

\begin{proof}
In \cite{r18} we proved this theorem using the Schwarz triangle function $s=s(\frac{2}{3},\frac{2}{3},\frac{2}{3},x)$ with an equilateral domain. Here we consider the Schwarz triangle function $s=s(\frac{4}{3},\frac{1}{3},\frac{1}{3},x)$ with a domain with isosceles symmetry and also the Schwarz triangle function $s=s(\frac{2}{3},\frac{1}{3},\frac{1}{2},x)$. 
For the Schwarz triangle function $s(\frac{4}{3},\frac{1}{3},\frac{1}{3},x)$, $y=P=-4w_1-w_2-w_3$ solves the $k=3$ generalised Chazy equation while $y=p=-2w_1-2w_2-2w_3$ solves $k=\frac{3}{2}$ generalised Chazy equation. 
For $y=P=-4w_1-w_2-w_3$ we obtain the triple
\begin{align*}
(P,Q,R)=&(-4w_1-w_2-w_3,-48(w_1-w_2)(w_1-w_3)-3(w_2-w_3)^2,\\
&3(2w_1-w_2-w_3)(32(w_1-w_3)(w_1-w_2)-(w_2-w_3)^2))
\end{align*}
and for
$y=p=-2w_1-2w_2-2w_3$ we obtain the triple
\begin{align*}
(p,q,r)=&(-2w_1-2w_2-2w_3,-60(w_1-w_3)(w_1-w_2),\\
&~60(w_1-w_3)(w_1-w_2)(2w_1-w_2-w_3)).
\end{align*}
Now we take $\mu=w_3-w_1$, $\la=w_2-w_3$, with $w_1-w_2=-\mu-\la$ and solve for $p$, $q$, $r$ in terms of $P$, $Q$, $R$.
We find
\begin{align*}
Q&=-48(-\mu-\la)(-\mu)-3\la^2=-48\mu^2-48\mu\la-3\la^2,\\
R&=-3(2\mu+\la)(-\frac{2}{3}Q-3\la^2).
\end{align*}
Let $w=2\mu+\la$. Then we have 
\begin{align*}
Q&=-12w^2+9\la^2,\\
R&=3wQ+12w^3.
\end{align*}
Consequently, we find $w^3+\frac{1}{4}Qw-\frac{1}{12}R=0$ and
\begin{align*}
(p,q,r)&=\left(P-2\mu-\la,-60(-\mu)(-\mu-\la),-60(-\mu)(-\mu-\la)(-2\mu-\la)\right)\\
&=\left(P-w,\frac{5}{3}Q+5w^2,-(\frac{5}{3}Qw+5w^3)w\right)\\
&=\left(P-w,\frac{5}{3}Q+5w^2,\frac{5}{12}R+\frac{5}{12}Qw\right).
\end{align*}

For the Schwarz triangle function $s(\frac{2}{3},\frac{1}{3},\frac{1}{2},x)$, $y=P=-2w_1-w_2-3w_3$ solves the $k=3$ generalised Chazy equation while $y=p=-w_1-2w_2-3w_3$ solves $k=\frac{3}{2}$ generalised Chazy equation. For  $y=P=-2w_1-w_2-3w_3$ we obtain the triple
\[
(P,Q,R)=(-2w_1-w_2-3w_3,-3(w_1-w_2)((w_1-w_2)+3(w_1-w_3)),3(w_1-w_2)^2(8(w_1-w_3)+w_2-w_3)).
\]
For $y=p=-w_1-2w_2-3w_3$ we obtain the triple
\begin{align*}
(p,q,r)=(-w_1-2w_2-3w_3,-15(w_1-w_3)(w_1-w_2),15(w_1-w_2)^2(w_1-w_3)).
\end{align*}
Now we take $\mu=w_3-w_1$, $\nu=w_1-w_2$ and solve for $p$, $q$, $r$ in terms of $P$, $Q$, $R$.
We find
\begin{align*}
Q=&-3\nu(-3\mu+\nu)=9\mu\nu-3\nu^2,\\
R=&3\nu^2(-\nu-9\mu)=-3\nu^3-27\mu\nu^2
=-3\nu^3-3(Q+3\nu^2)\nu
=-12\nu^3-3Q\nu.
\end{align*}
This gives $\nu^3+\frac{1}{4}Q\nu+\frac{1}{12}R=0$
and
\begin{align*}
(p,q,r)=&\left(P+\nu,15\mu\nu,-15\nu^2\mu\right)\\
=&\left(P+\nu,\frac{5}{3}Q+5\nu^2,\frac{5}{9}(R+3\nu^3)\right)\\
=&\left(P+\nu,\frac{5}{3}Q+5\nu^2,\frac{5}{9}(\frac{3}{4}R-\frac{3}{4}Q\nu)\right)\\
=&\left(P+\nu,\frac{5}{3}Q+5\nu^2,\frac{5}{12}R-\frac{5}{12}Q\nu\right).
\end{align*}
Substituting $w=-\nu$ gives  $w^3+\frac{1}{4}Qw-\frac{1}{12}R=0$
and
\begin{align*}
(p,q,r)=&\left(P-w,\frac{5}{3}Q+5w^2,\frac{5}{12}R+\frac{5}{12}Qw\right),
\end{align*}
which we have already deduced earlier. 
Conversely, starting from the triple $(P,Q,R)$ a solution to the system of differential equations associated to the $k=\frac{3}{2}$ generalised Chazy equation, we obtain straightforwardly that the triple
\begin{align*}
\left(P+\frac{R}{Q},\frac{3}{5}Q-3\frac{R^2}{Q^2},\frac{9}{5}R+3\frac{R^3}{Q^3}\right)
\end{align*}
satisfies the first order system associated to the $k=3$ generalised Chazy equation.

\end{proof}

\section{Automorphism of solutions to $k=4$ generalised Chazy equation and transformation to solutions of $k=3$ generalised Chazy equation}\label{aut4}

In this section we find an automorphism of the solution to the generalised Chazy equation with parameter $k=4$ and a transformation to the solution of the generalised Chazy equation with parameter $k=3$.
\begin{theorem}\label{thm4}
Suppose $(P,Q,R)$ satisfies the non-linear system of differential equations associated to the generalised Chazy equation with parameter $k=4$, i.e. we have
\begin{align*}
\frac{\der}{\der x}P&=\frac{1}{6}(P^2-Q),\nonumber\\
\frac{\der}{\der x}Q&=\frac{2}{3}(PQ-R),\\
\frac{\der}{\der x}R&=PR+\frac{4}{5}Q^2.\nonumber
\end{align*}
Then the following holds. The triple
\begin{align*}
\left(P-\la,-\frac{5}{3}Q-5\la^2-\frac{10}{3}\nu^2,-\frac{5}{3}Q\nu-\frac{20}{27}R-7Q\la-8\frac{Q^2}{R}\la^2-\frac{8}{135}\frac{Q^3}{R}\right)
 \end{align*}
also satisfies the non-linear system of differential equations associated to the generalised Chazy equation with parameter $k=4$, where $\nu$ is a root of the cubic equation
\[
\nu^3+\frac{3}{5}Q\nu-\frac{1}{5}R=0
\]
and $\la=\frac{1}{3}\nu+\frac{4}{15\nu}Q$ satisfies
\[
\la^3-\frac{4}{5}\frac{Q^2}{R}\la^2-\frac{Q}{5}\la-\frac{4}{675}\frac{Q^3}{R}-\frac{1}{135}R=0.
\]
Furthermore, let 
\[
\chi=-\frac{R}{2Q}\pm\frac{1}{2}\sqrt{\frac{R^2}{Q^2}+\frac{4}{5}Q}.
\]
Then the triple
\[
\left(P+\chi, -\frac{3}{5}Q-3\chi^2,-\frac{3}{5}R-\frac{9}{5}Q\chi-3\chi^3\right)
\]
satisfies the first order system associated to the $k=3$ generalised Chazy equation. 
\end{theorem}
\begin{proof}
To prove this theorem we first consider the Schwarz triangle function $s(\frac{1}{4},\frac{1}{2},\frac{1}{2},x)$. We find that $y=P=-w_1-2w_2-3w_3$ and $y=p=-w_1-3w_2-2w_3$ both solve the $k=4$ generalised Chazy equation. 
For $y=P=-w_1-2w_2-3w_3$, using (\ref{wde}) with $(\al,\beta,\ga)=(\frac{1}{4},\frac{1}{2},\frac{1}{2})$ and $Q=P^2-6P'$, $R=PQ-\frac{3}{2}Q'$, we find that
\begin{align*}
Q=&-\frac{5}{4}(w_1-w_2)((w_1-w_2)+3(w_3-w_2)),\\
R=&\frac{5}{4}(w_1-w_2)^2((w_1-w_3)+8(w_2-w_3)),
\end{align*}
while for $y=p=-w_1-3w_2-2w_3$, with $q=p^2-6p'$ and $r=pq-\frac{3}{2}q'$ we have
\begin{align*}
p=&-w_1-3w_2-2w_3,\\
q=&-\frac{5}{4}(w_1-w_3)((w_1-w_3)+3(w_2-w_3)),\\
r=&\frac{5}{4}(w_1-w_3)^2((w_1-w_2)+8(w_3-w_2)).
\end{align*}
Let $\la=w_2-w_3$, $\nu=w_1-w_2$ and $w_1-w_3=\nu+\la$.
Then the formulas for $Q$ and $R$ give
\begin{align*}
Q=&-\frac{5}{4}\nu(\nu-3\la)=-\frac{5}{4}\nu^2+\frac{15}{4}\nu\la,\\
R=&\frac{5}{4}\nu^2(9\la+\nu)=\frac{5}{4}\nu^3+\frac{45}{4}\nu^2\la
=\frac{5}{4}\nu^3+3Q\nu+\frac{15}{4}\nu^3
=5\nu^3+3Q\nu,
\end{align*}
from which we obtain
$5\nu^3+3Q\nu-R=0$ and $\la=\frac{1}{3}\nu+\frac{4}{15}\frac{Q}{\nu}$ where $\la$ satisfies the cubic equation
\[
\la^3-\frac{4}{5}\frac{Q^2}{R}\la^2-\frac{Q}{5}\la-\frac{4}{675}\frac{Q^3}{R}-\frac{1}{135}R=0.
\]
We also find
\begin{align*}
\nu^2\la&=\frac{1}{15}(Q\nu+R),\\
\nu\la^2&=\frac{4}{15}Q\la+\frac{1}{45}Q\nu+\frac{1}{45}R.
\end{align*}
This allows us to express the triple $(p,q,r)$ satisfying the first order system associated to the $k=4$ generalised Chazy equation, with
\begin{align*}
p=&P-\la,\\
q=&-\frac{5}{4}(\nu+\la)(\nu+4\la)
=-\frac{5}{4}(\nu^2+5\la\nu+4\la^2)
=-\frac{5}{4}(\nu^2+\frac{4}{3}(Q+\frac{5}{4}\nu^2)+4\la^2)\\
=&-\frac{5}{4}(\frac{8}{3}\nu^2+\frac{4}{3}Q+4\la^2)=-\frac{10}{3}\nu^2-\frac{5}{3}Q-5\la^2,\\
r=&\frac{5}{4}(\nu+\la)^2(\nu-8\la)=\frac{5}{4}(\nu^3-6\la\nu^2-15\la^2\nu-8\la^3)\\
=&-\frac{5}{3}Q\nu-\frac{20}{27}R-7Q\la-8\frac{Q^2}{R}\la^2-\frac{8}{135}\frac{Q^3}{R}.
\end{align*}
To obtain the transformation to the solution of the $k=3$ generalised Chazy equation, once again consider the Schwarz triangle function $s=s(\frac{1}{3},\frac{1}{3},\frac{1}{2},x)$ as in Section \ref{aut3}. The combination $y=P=-2w_1-2w_2-2w_3$ solves the $k=4$ generalised Chazy equation and $y=p=-w_1-2w_2-3w_3$ and $y=\tilde p=-2w_1-w_2-3w_3$ both solve the $k=3$ generalised Chazy equation. We have
\begin{align*}
(P,Q,R)=(-2w_1-2w_2,-2w_3,-5(w_1-w_3)(w_2-w_3),-5(w_1-w_3)(w_2-w_3)(w_1+w_2-2w_3)).
\end{align*}
We also obtain the triples
 \[
(p,q,r)=(-2w_1-w_2-3w_3,3(w_2-w_3)(w_1-w_2),3(w_1-w_2)^2(w_2-w_3))
\]
and
\[
(\tilde p,\tilde q,\tilde r)=(-w_1-2w_2-3w_3,-3(w_1-w_3)(w_1-w_2),3(w_1-w_2)^2(w_1-w_3))
\]
for $p=-2w_1-w_2-3w_3$ and $\tilde p=-w_1-2w_2-3w_3$.
Let $\la=w_2-w_3$ and $w_1-w_3=-\mu$, with $w_1-w_2=-\mu-\la$. 
Then 
\begin{align*}
\frac{R}{Q}&=w_1+w_2-2w_3=\la-\mu,\\
Q&=-5(-\mu)\la=5\la\mu. 
\end{align*}
Solving the quadratic equation
\[
\chi^2+\frac{R}{Q}\chi-\frac{Q}{5}=0
\]
for $\chi=\mu$ and $\chi=-\la$ gives
\[
\chi=-\frac{R}{2Q}\pm\frac{1}{2}\sqrt{\frac{R^2}{Q^2}+\frac{4}{5}Q}.
\]
Therefore for either triples $(p,q,r)$ or $(\tilde p, \tilde q, \tilde r)$ we get
\[
\left(P+\chi, -\frac{3}{5}Q-3\chi^2,-\frac{3}{5}R-\frac{9}{5}Q\chi-3\chi^3\right)
\]
as solutions to the first order system associated to the $k=3$ generalised Chazy equation. 

\end{proof}

\section{Automorphism of solutions to $k=9$ generalised Chazy equation and transformation to solution of the $k=18$ generalised Chazy equation}\label{aut9}

This section concerns the automorphism of solutions to the $k=9$ generalised Chazy equation and its transformation to the solution of the $k=18$ generalised Chazy equation.
 
\begin{theorem}\label{thm9}
Suppose $(P,Q,R)$ satisfies the non-linear system of differential equations associated to the generalised Chazy equation with parameter $k=9$, i.e. we have
\begin{align*}
\frac{\der}{\der x}P&=\frac{1}{6}(P^2-Q),\nonumber\\
\frac{\der}{\der x}Q&=\frac{2}{3}(PQ-R),\\
\frac{\der}{\der x}R&=PR-\frac{9}{5}Q^2.\nonumber
\end{align*}
Then the following holds. The triple
\begin{align*}
\left(P-\la,-\frac{5}{4}Q+\frac{5}{36}\mu^2+5\la^2,\frac{5}{9}\mu^3-15\la^3+\frac{349}{54}Q\la+\frac{181}{216}Q\mu-\frac{181}{216}R\right)
 \end{align*}
also satisfies the non-linear system of differential equations associated to the generalised Chazy equation with parameter $k=9$, where $\mu$ is a root of the quartic equation
\begin{align*}
\mu^4-\frac{6}{5} Q\mu^2-\frac{8}{15}R\mu-\frac{3}{25}Q^2=0,
\end{align*}
and $\la=\frac{1}{8}\mu-\frac{9}{40}\frac{Q}{\mu}$ satisfies
\[
\la^4-\frac{R}{Q}\la^3+\frac{3}{5}Q\la^2-\frac{R}{15}\la+\frac{1}{960}\frac{R^2}{Q}-\frac{3}{1600}Q^2=0.
\]
For the same $\mu$ and $\la$ the triple
\begin{align*}
\left(P-\mu,-\frac{4}{5}Q+4\mu^2,\frac{4}{15}R+8\mu^3-\frac{28}{15}Q\mu-\frac{16}{15}Q\la\right)
 \end{align*}
satisfies the non-linear system of differential equations associated to the generalised Chazy equation with parameter $k=18$.
Conversely, let $(P,Q,R)$ be a solution to the $k=18$ generalised Chazy equation, i.e. we have
\begin{align*}
\frac{\der}{\der x}P&=\frac{1}{6}(P^2-Q),\nonumber\\
\frac{\der}{\der x}Q&=\frac{2}{3}(PQ-R),\\
\frac{\der}{\der x}R&=PR-\frac{9}{8}Q^2.\nonumber
\end{align*}
Let 
\[
\chi=-\frac{R}{2Q}\mp\sqrt{\frac{R^2}{4Q^2}-\frac{9}{32}Q}.
\]
Then 
\[
\left(P+\chi,-\frac{5}{32}Q+5\chi^2,-\frac{5}{4}R-\frac{55}{8}Q\chi+15\chi^3\right)
\]
solves the first order system associated to the $k=9$ generalised Chazy equation.
\end{theorem}
\begin{proof}
Here we consider the solution to the generalised Chazy equation with parameter $k=9$ given by the Schwarz triangle function $s(\frac{1}{9},\frac{1}{3},\frac{1}{3},x)$. We find $y=P=-w_1-2w_2-3w_3$ and $y=p=-w_1-3w_2-2w_3$ both solve the $k=9$ generalised Chazy equation, and that $y=\tilde p=-2w_1-2w_2-2w_3$ solves the $k=18$ generalised Chazy equation.  
For $y=P=-w_1-2w_2-3w_3$, we find from (\ref{wde}) with $(\al,\beta,\ga)=(\frac{1}{9},\frac{1}{3},\frac{1}{3})$
that we have
\begin{align*}
P=&-w_1-2w_2-3w_3,\\
Q=&\frac{5}{9}(w_1-w_3)((w_1-w_3)+8(w_2-w_3)),\\
R=&-\frac{5}{9}(w_1-w_3)((w_1-w_2)^2-9(w_2-w_3)^2+18(w_1-w_3)(w_3-w_2)),
\end{align*}
while for $y=p=-w_1-3w_2-2w_3$, we find
\begin{align*}
p=&-w_1-3w_2-2w_3,\\
q=&\frac{5}{9}(w_1-w_2)((w_1-w_2)+8(w_3-w_2)),\\
r=&-\frac{5}{9}(w_1-w_2)((w_1-w_3)^2-9(w_2-w_3)^2+18(w_1-w_2)(w_2-w_3)).
\end{align*}
For the solution to the $k=18$ generalised Chazy equation given by $y=\tilde p=-2w_1-2w_2-2w_3$, we obtain
\begin{align*}
\tilde p=&-2w_1-2w_2-2w_3,\\
\tilde q=&\frac{32}{9}(w_1-w_2)((w_1-w_3),\\
\tilde r=&-\frac{32}{9}(w_1-w_2)(w_1-w_3)(2w_1-w_2-w_3).
\end{align*}
Let $\la=w_2-w_3$, $w_1-w_3=-\mu$ and $w_1-w_2=-\mu-\la$.
Then 
\begin{align*}
Q=&-\frac{5}{9}\mu(-\mu+8\la)=\frac{5}{9}\mu(\mu-8\la),\\
R=&\frac{5}{9}\mu((\mu+\la)^2-9\la^2+18\la\mu)=\frac{5}{9}\mu(\mu^2-8\la^2+20\la\mu)=Q\mu-\frac{5}{9}(8\la^2-28\la\mu).
\end{align*}

Eliminating $\mu$ in the equations for $Q$ and $R$ for $\la$, we find that $\la$ satisfies the quartic equation
\[
\la^4-\frac{R}{Q}\la^3+\frac{3}{5}Q\la^2-\frac{R}{15}\la+\frac{1}{960}\frac{R^2}{Q}-\frac{3}{1600}Q^2=0.
\]
Furthermore, we have
\[
\la=\frac{1}{8}\mu-\frac{9}{40}\frac{Q}{\mu}
\]
where $\mu$ satisfies the quartic equation
\[
\mu^4-\frac{6}{5}Q\mu^2-\frac{8}{15}R\mu-\frac{3}{25}Q^2=0.
\]

From the equations for $Q$ and $R$ we also have
\begin{align*}
Q\la&=\frac{5}{9}\mu^2\la-\frac{40}{9}\mu\la^2,\\
Q\mu-R&=\frac{40}{9}\mu\la^2-\frac{140}{9}\la\mu^2.
\end{align*}
Now from the formulas for $p$, $q$ and $r$ we have
\begin{align*}
p=&P-\la,\\
q=&\frac{5}{9}(-\mu-\la)((-\mu-\la)-8\la)=\frac{5}{9}(\mu+\la)(9\la+\mu)=\frac{5}{9}(\mu^2+9\la^2+10\mu\la),\\
r=&-\frac{5}{9}(-\mu-\la)(\mu^2-9\la^2+18(-\mu-\la)\la)\\
=&\frac{5}{9}(\mu+\la)(\mu^2-27\la^2-18\mu\la)\\
=&\frac{5}{9}(\mu^3-45\la^2\mu-17\mu^2\la-27\la^3).
\end{align*}
We therefore obtain
\begin{align*}
p=&P-\la,\\
q=&-\frac{5}{4}Q+\frac{5}{36}\mu^2+5\la^2,\\
r=&\frac{5}{9}(\mu^3-27\la^3+\frac{17}{15}(Q\la+Q\mu-R)+\frac{3}{8}(28Q\la+Q\mu-R))\\
=&\frac{5}{9}(\mu^3-27\la^3+\frac{349}{30}Q\la+\frac{181}{120}Q\mu-\frac{181}{120}R)\\
=&\frac{5}{9}\mu^3-15\la^3+\frac{349}{54}Q\la+\frac{181}{216}Q\mu-\frac{181}{216}R,
\end{align*}
and this gives the automorphism
\[
(p,q,r)=\left(P-\la,-\frac{5}{4}Q+\frac{5}{36}\mu^2+5\la^2,\frac{5}{9}\mu^3-15\la^3+\frac{349}{54}Q\la+\frac{181}{216}Q\mu-\frac{181}{216}R\right)
\]
to another solution of the $k=9$ generalised Chazy equation. Similarly, we have
\begin{align*}
\tilde p=&P-(-\mu-\la)-\la=P+\mu,\\
\tilde q=&\frac{32}{9}(-\mu-\la)(-\mu)=\frac{32}{9}(\mu^2+\la\mu)=-\frac{4}{5}Q+4\mu^2,\\
\tilde r=&-\frac{32}{9}(-\mu-\la)(-\mu)(2(-\mu-\la)+\la)=\frac{32}{9}\mu(\mu+\la)(2\mu+\la)\\
=&\frac{4}{15}R+8\mu^3-\frac{28}{15}Q\mu-\frac{16}{15}Q\la
\end{align*}
and we obtain 
\begin{align*}
(\tilde p, \tilde q, \tilde r)=\left(P-\mu,-\frac{4}{5}Q+4\mu^2,\frac{4}{15}R+8\mu^3-\frac{28}{15}Q\mu-\frac{16}{15}Q\la\right)
\end{align*}
as a solution of the first order system associated to the $k=18$ generalised Chazy equation.

Conversely, given $(P,Q,R)$ a solution to the system of differential equations associated to the $k=18$ generalised Chazy equation, we find
\begin{align*}
\nu-\mu&=2w_1-w_2-w_3=\frac{R}{Q}\\
-\mu\nu&=(w_1-w_2)(w_1-w_3)=\frac{9}{32}Q,
\end{align*}
where $\nu=w_1-w_2$ and $\mu=w_3-w_1$. 
Solving these two equations for $\chi=-\nu$ or $\chi=\mu$ gives us
\[
\chi=-\frac{R}{2Q}\mp\sqrt{\frac{R^2}{4Q^2}-\frac{9}{32}Q}
\]
and we have the triple
\[
\left(P+\chi,-\frac{5}{32}Q+5\chi^2,-\frac{5}{4}R-\frac{55}{8}Q\chi+15\chi^3\right)
\]
satisfying the first order system of dfferential equations associated to the $k=9$ generalised Chazy equation.
\end{proof}

The next step of the project would be to investigate composition of the automorphism of the solutions and transformations of solutions between various different Chazy parameters. For example, we can express solutions to the $k=\frac{3}{2}$ generalised Chazy equation in terms of solutions to the $k=\frac{2}{3}$ solutions via composing through the transformation to $k=2$ and $k=3$ solutions. Also, the solutions to the $k=2$, $k=3$ and $k=4$ generalised Chazy equation have dihedral, tetrahedral and octahedral symmetry respectively and the transformations between these solutions establishes a link between these three symmetries.

\end{document}